\documentclass[12pt,twoside,a4paper,reqno]{amsart}

\usepackage[a4paper,margin=1.22in,footskip=0.85in]{geometry}

\usepackage{amsmath,amsthm,amsfonts,amssymb}
\usepackage{amsthm}
\usepackage{graphicx}
\usepackage{hyperref}
\usepackage{color}
\usepackage[T1]{fontenc}
\usepackage[utf8]{inputenc}
\usepackage{xy}
\usepackage{graphics}
\usepackage{mathrsfs}
\usepackage{comment}
\usepackage{xcolor}
\usepackage{comment}
\xyoption{all}

\makeatletter
\newtheorem*{rep@theorem}{\rep@title}
\newcommand{\newreptheorem}[2]{%
\newenvironment{rep#1}[1]{%
 \def\rep@title{#2 \ref{##1}}%
 \begin{rep@theorem}}%
 {\end{rep@theorem}}}
\makeatother

\numberwithin{equation}{section}
\theoremstyle{plain}
\newtheorem{theorem}{Theorem}[section]
\newreptheorem{theorem}{Theorem}
\newtheorem{proposition}[theorem]{Proposition}
\newtheorem{lemma}[theorem]{Lemma}
\newtheorem{corollary}[theorem]{Corollary}
\newreptheorem{corollary}{Corollary}

\theoremstyle{definition}
\newtheorem{definition}[theorem]{Definition}

\newtheorem{remark}[theorem]{Remark}

\newtheorem{example}[theorem]{Example}

\usepackage{tikz-cd}
\title{Stability of syzygy bundles over certain curves of compact type}
\author[Amit Kumar Singh]{Amit Kumar Singh}
\address{Department of Mathematics, SRM University AP, Amaravati 522240, Andhra Pradesh, India}
\email{amitks.math@gmail.com}
\begin{document}
\maketitle
\begin{abstract}Let $E$ be a vector bundle on a curve $C$ of compact type and $V \subseteq \mathrm{H}^0(C, E)$ be a linear subspace that generates $E$. In this note, we study the stability of the syzygy bundle $M_{E,V}$ associated to $(E, V )$ over certain curves $C$ of compact type.
\end{abstract}
\section{Introduction}
All the curves considered throughout this note are projective curves and defined over the complex field. Let $E$ be a rank $r$ vector bundle on a curve $C$ of compact type and $V$ be an $ l$-dimensional linear subspace of $\mathrm{H}^0(C, E)$ with $l>r$ that generates $E$ (that means the evaluation map $V \otimes \mathcal{O}_C \xrightarrow{\mathrm{ev}} E$ is onto). Such a pair $(E,V)$ is called a {\em generated pair} on $C$. The {\em syzygy bundle} $M_{E, V}$ on $C$ of the generated pair $(E, V)$ is defined by the following exact sequence:
\begin{eqnarray}\label{exact seq defining seq for kernel}
0 \rightarrow M_{E, V} \rightarrow V\otimes \mathcal{O}_C \xrightarrow{\mathrm{ev}} E \rightarrow 0.
\end{eqnarray} 
In case if $V = \mathrm{H}^0(C, E)$, then we denote the syzygy bundle by $M_E$ instead of $M_{E, V}$. For a generated pair $(E,V)$ on $C$, we have the following map:
\begin{eqnarray*}
\varphi_{E,V} : C \rightarrow \mathrm{Gr}(E, V)  := \mathrm{Gr}
\end{eqnarray*}
defined by setting $\varphi_{E,V}(x) = E_x$ and then we have a canonical exact sequence:
\begin{eqnarray}\label{exact seq on universival on grassmannian}
0 \rightarrow S \rightarrow V\otimes \mathcal{O}_\mathrm{Gr} \rightarrow Q \rightarrow 0
\end{eqnarray}
with $\varphi_{E,V}^*(S) \cong M_{E, V}$. In particular, if the vector bundle $E$ is a line bundle $L$ (say), then \eqref{exact seq on universival on grassmannian} is nothing but the Euler exact sequence
\begin{eqnarray*}
0 \rightarrow \Omega_\mathbb{P}(1) \rightarrow V \otimes \mathcal{O}_\mathbb{P} \rightarrow \mathcal{O}_\mathbb{P}(1) \rightarrow 0
\end{eqnarray*}
on $ \mathbb{P}(V^*) = \mathbb{P}$. Moreover $M_{L, V} \cong \varphi^*_{L,V}(\Omega_{\mathbb{P}}(1))$, where $\Omega_{\mathbb{P}}(1)$ denotes the twisted cotangent bundle on $\mathbb{P}$.

One of the well-explored questions in algebraic geometry is to study the stability of the syzygy bundles as it has a deep connection with higher rank Brill Noether Theory and theta-divisor (see \cite{Kapil-thesis, Hein98, Beauville2006, Mistretta08, Popa99}).

Over a smooth curve $C$ of genus $g \ge 2$, Ein and Lazarsfeld proved that $M_{L}$ is (semi)stable provided $\deg L (\ge) > 2g$ \cite[Proposition 3.1]{EL92}. Later, this result was generalised for (semi)stable vector bundles $E$ of rank $r$ on $C$ by Butler. More precisely, he proved that $M_E$ is (semi)stable provided $\deg E (\ge) > 2 g r$ \cite[Theorem 1.2]{Butler94}. Much later, Camere sharpened the bound on the degree for line bundle by using the Clifford index \cite[Theorem 1.3]{Camere2008}.

A considerable amount of work has been done to tackle the question about the stability of $M_{E, V}$ if $V \neq \mathrm{H}^0(C, E)$ (see, for example \cite{Butler-conj, Usha08, Usha15, Brambila-Newstead2019}). Nevertheless, it has not yet reached a complete answer of the question. In 2015, Bhosle, Brambila-Paz and Newstead proved $M_{L,V}$ is semistable for a general linear system $(L, V)$ over a general curve $C$ with genus $g \ge 1$ by using wall crossing formulae for coherent systems \cite[Theorem 5.1]{Usha15}.

Though there is a rich literature on the study of the stability of syzygy bundles over a smooth curve, still not much is known for the case of a singular curve. On certain reducible nodal curves with ordinary nodes, the stability of $M_{E,V}$ was studied in \cite{Brivio-Favale2020, SSA, SSA2025}. Surprisingly, the obtained results go in the opposite direction with to what is known for the smooth curve cases.

In this note, we prove that the syzygy bundle $M_{E,V}$ is strongly unstable on a certain curve of compact type (say comb-like curve) assuming some conditions (see Theorem \ref{thm on strongly unstabe for any rank}). 
\section{Syzygy bundles on a curve of compact type}
Let C be a curve of compact type having $n \ge 2$ smooth components $C_1, \dots , C_n$ with genus $g_i$.  On such a curve, we have a canonical short exact sequence:
\begin{equation}\label{canonical exact sequence}
    0 \rightarrow \mathcal{O}_C \rightarrow \bigoplus_{i=1}^{n
    }\mathcal{O}_{C_i} \rightarrow \mathcal{T} \rightarrow 0,
 \end{equation}
 where the skyscraper sheaf $\mathcal{T}$ is supported only at the nodal point(s). And at each nodal point $p$, $\mathcal{T}_{p}$ is a one dimensional complex vector space and that $\dim \mathrm{H}^0(\mathcal{T}) = n-1$.
 From the exact sequence \ref{canonical exact sequence}, we have 
\begin{eqnarray*}
    \chi(\mathcal{O}_C) & = & \sum_{i=1}^{n} \chi(\mathcal{O}_{C_i}) - \chi(\mathcal{T}) \nonumber \\
    & = & n - \sum_{i=1}^{n} g_i - (n-1)
    \nonumber \\
    & = & 1 - \sum_{i=1}^n g_i
\end{eqnarray*}
So if we let $p_a(C) = 1 - \chi(\mathcal{O}_C)$ be the arithmetic genus of $C$, then 
\begin{eqnarray}\label{exact seq calculation of arithmetic genus}
p_a(C) = \sum\limits_{i=1}^{n} g_i.
\end{eqnarray}

Now, let $(E, V)$ be a generated pair on $C$. For each $i$ (we denote the restriction $E_{\vert_{C_i}}$ of $E$ to $C_i$ by $E_i$), we have a natural restriction map
\[
\rho_i : \mathrm{H}^0 (C, E) \rightarrow \mathrm{H}^0 (C_i, E_i)  
\]
defined by $\rho_i(s) = s_{\vert_{C_i}}$. The image $\rho_i(V)$ of $V$ under $\rho_i$ is denoted by $V_i$. The pair $(E_i, V_i)$ is said to be of type $(r, d_i, l_i)$ where $r = \mathrm{rk}(E_i)$, $d_i = \deg(E_i)$ and $l_i = \dim(V_i)$.  Then we have the following diagram:
\begin{eqnarray}\label{diagram 1}
\begin{tikzcd}
  & 0 \arrow[d, ]  & 0 \arrow[d, ]  & 0 \arrow[d, ]  \\
   0 \arrow[r] & {\ker\big(\rho_{i_{|_V}} \big)} \otimes \mathcal{O}_{C_i} \arrow[equal]{d} \arrow[r, ] & {M_{E,V}}_{\vert_{C_i}} \arrow[d, ] \arrow[r] & M_{E_i, V_i} \arrow[d, ] \arrow[r] & 0 \\
    0 \arrow[r] & {\ker \big(\rho_{i_{|_V}} \big)} \otimes \mathcal{O}_{C_i} \arrow[d, ] \arrow[r, ] & V\otimes \mathcal{O}_{C_i} \arrow[d, ] \arrow[r,] & V_i \otimes \mathcal{O}_{C_i} \arrow[d, ] \arrow[r] & 0 \\
     & 0  \arrow[r, ] & E_i \arrow[d, ] \arrow[equal]{r} & E_i \arrow[d, ] \arrow[r] & 0\\
  & &0  & 0   
\end{tikzcd}
\end{eqnarray}
Suppose that $E$ be a rank $r$ vector bundle on $C$. After tensoring in \eqref{canonical exact sequence} with $E$,
we have 
\begin{eqnarray}\label{exact seq tensor with E}
0 \rightarrow E \rightarrow \bigoplus\limits_{i = 1}^n {E_i} \rightarrow \tau_E \rightarrow 0.
\end{eqnarray}
Then by using the additive property of the Euler characteristic, we obtain
\begin{eqnarray*}
\begin{split}
\chi(E) & = \sum\limits_{i=1}^n \chi(E_i) - r(n-1)\\
        & = \sum\limits_{i=1}^n d_i + r(1-p_a(C)).
\end{split}
\end{eqnarray*}

\begin{remark}
Let $(E,V)$ be a generated pair on $C$. Then from the defining exact sequence \eqref{exact seq defining seq for kernel} of $M_{E,V}$, the syzygy bundle $M_{E,V}$ is a subbundle of the trivial bundle $V \otimes \mathcal{O}_{C}$ and we have
\[
\begin{split}
\chi(M_{E,V}) & = \chi(V \otimes \mathcal{O}_{C}) - \chi(E) \\
           & = l(1-p_a(C)) - \sum\limits_{i=1}^n d_i - r(1-p_a(C)) \\
           & = (l-r)(1-p_a(C)) - \sum\limits_{i=1}^n d_i.
\end{split}
\]
Therefore,
\begin{eqnarray*}
\chi(M_{E,V}) = (l-r)(1-p_a(C)) -  d,
\end{eqnarray*}
where $d = \sum\limits_{i=1}^n d_i$. Moreover, $\chi(M_{E, V}) < 0$.
\end{remark}

\begin{proposition}
Let $(E,V)$ be a generated pair on $C$. Then
\begin{enumerate}
\item[(a)] for each $i$, the pair $(E_i, V_i)$ is generated pair on $C_i$.
\item[(b)] $\dim\mathrm{H}^0(C,E) = \sum\limits_{i=1}^n \dim\mathrm{H}^0(C, E_i) + (n-1)r$.
\end{enumerate}
\begin{proof}
The first result (a) follows from the above diagram \ref{diagram 1}.

\noindent
(b) Consider the exact sequence \eqref{exact seq tensor with E} and, then by passing the cohomology, we have
\begin{eqnarray*}
0 \rightarrow \mathrm{H}^0(C, E) \rightarrow \bigoplus\limits_{i = 1}^n \mathrm{H}^0(C _i, {E_i}) \rightarrow \mathrm{H}^0(\tau_E) \rightarrow 0.
\end{eqnarray*} We conclude the result by calculating the dimension of the vector spaces in the above exact sequence.

\end{proof}
\end{proposition}

 Let $F$ be a coherent sheaf of $\mathcal{O}_C$-modules. We call $F$ a pure sheaf of dimension one if for every proper $\mathcal{O}_C$-submodule  $G \subset F$ and $G \neq 0$, the dimension of support of $G$ is equal to one. Vector bundles on $C$ are examples of pure sheaves of dimension one. Suppose $F$ is a pure sheaf of dimension one on $C$. Let $F_i = \frac{F_{|_{C_i}}}{\text{Torsion}(F_{|_{C_i}})}$ for each $i$, where $\text{Torsion}(F_{|_{C_i}})$ is the torsion subsheaf of $F_{|_{C_i}}$. Then $F_i$, if non-zero, is torsion-free and hence locally free on $C_i$. Let $r_i$ denote the rank of $F_i$. Also, let $d_i$ denote the degree of $F_i$ for each $i$. Then we call the n-tuples $(r_1,\dots,r_n)$ and $(d_1,\dots,d_n)$ respectively the \textit{multirank} and \textit{multidegree} of $F$.
\begin{definition}\label{polarization}
 Let $w = (w_1,\dots,w_n)$ be an $n$-tuple with $ w_i \in (0,\;1) \cap \mathbb{Q}$ for each $i$ and $\sum\limits_{i=1}^n w_i =1$. We call such an $n$-tuple a polarization on $C$.
\end{definition}

\begin{definition}\label{polarized slope}
 Suppose $F$ is a pure sheaf of dimension one on $C$ of multirank $(r_1,\dots,r_n)$. Then the slope of $F$ with respect to a polarization $w$, denoted by $\mu_w(F)$, is defined by $\mu_w(F) = \frac{\chi(F)}{\sum\limits_{i=1}^n w_i \; r_i}$.
\end{definition}

\begin{definition}
 Let $E$ be a vector bundle defined on $C$. Then $E$ is said to be $w$-semistable (resp. $w$-stable) if for any proper subsheaf $F \subset E$ one has $\mu_w(F) \leq \mu_w(E)$ (resp. $\mu_w(F) < \mu_w(E)$). If $E$ is not $w$-semistable for any polarization $w$ on $C$, we call it a strongly unstable bundle.
\end{definition}

\noindent
At the end of this section, we prove the following result for a curve of compact type. 

\begin{proposition}
Let $(E,V)$ be a generated pair on $C$. If $\ker({\rho_{i_{\vert_V}}}) \neq 0$, for some $i$, then ${M_{E,V}}_{\vert_{C_i}}$ is unstable.

\begin{proof}
Being $\deg (\ker ({\rho_i}_{\vert_V}) \otimes \mathcal{O}_{C_i}) =0$, the slope $\mu(\ker({\rho_{i_{\vert_V}}}) \otimes \mathcal{O}_{C_i}) = 0$. Now, then

\[
\mu(\ker({\rho_{i_{\vert_V}}}) \otimes \mathcal{O}_{C_i}) = 0 > \frac{-d_i}{l-r} = \mu\big({M_{E,V}}_{\vert_{C_i}}\big).
\]
Thus, we get a proper subbundle $\ker({\rho_{i_{\vert_V}}}) \otimes \mathcal{O}_{C_i}$ of ${M_{E,V}}_{\vert_{C_i}}$ satisfying
$\mu(\ker({\rho_{i_{\vert_V}}}) \otimes \mathcal{O}_{C_i}) > \mu\big({M_{E,V}}_{\vert_{C_i}} \big)$. This concludes the prove that $\mu\big({M_{E,V}}_{\vert_{C_i}} \big)$ is unstable.
\end{proof}
\end{proposition}

\section{Stability of syzygy bundles}
In this section, we prove the main result of this article. Moreover, we learn the strongly unstability of the syzygy bundle $M_{E,V}$ associated with a generated pair $(E, V)$ over a certain curve of compact type (say, comb-like curve) by assuming some restrictions on the invariants of the pair $(E,V)$.

Let $n\geq 2$ be an integer and let $C$ be a curve of compact type having $n$ smooth components $C_i$ of genus $g_i$ and $n-1$ nodes $p_i$ such that $C_i \cap C_j = \emptyset$ for $i \in \{1, \dots, n\}$, $j \notin \{i, n\}$ and $C_{i} \cap C_n = \{ p_i \}$, for $i \in \{1, \dots, n-1\}$. We call such a curve a comb-like curve with the fixed component $C_{n}$. On such a curve $C$, there is an exact sequence:
\begin{eqnarray}\label{exact sequence defining I_C_i}
0 \rightarrow \mathcal{I}_{C_i} \rightarrow \mathcal{O}_C \rightarrow \mathcal{O}_{C_i} \rightarrow 0,
\end{eqnarray}
where \[ \mathcal{I}_{C_i} \supseteq  (\bigoplus\limits_{\substack{j=1 \\ j \neq i}}^{n-1} \mathcal{O}_{C_j}(-p_j)) \oplus \mathcal{O}_{C_n}(-p_1 - \cdots - p_{n-1}) \quad \text{for} \; i \in \{ 1, \dots , n-1 \}\] and 
\[
\mathcal{I}_{C_n} = \bigoplus\limits_{i=1}^{n-1} \mathcal{O}_{C_i}(-p_i).
\]

\begin{remark}
Let $(E, V)$ be a generated pair on a comb-like curve $C$. Then, tensoring in \eqref{exact sequence defining I_C_i} by $E$, we have 
\begin{equation*} \label{exact sequence defining I_C_i tensor E}
 0 \rightarrow E \otimes \mathcal{I}_{C_i} \rightarrow E \rightarrow E_i \rightarrow 0.
\end{equation*}
Now, passing to the long exact sequence in the cohomology, we obtain
\begin{equation*}
0 \rightarrow \mathrm{H}^0(E \otimes \mathcal{I}_{C_i}) \rightarrow \mathrm{H}^0(E) \xrightarrow{\rho_i} H^0(E_i) \rightarrow \dots
\end{equation*}
This implies, in particular, 
\begin{eqnarray}\label{ker:N}
\text{ker}\big(\rho_{{n}_{|_V}} \big) = \bigoplus\limits_{i=1}^{n-1} (V \cap H^0(E_i(-p_i))),
\end{eqnarray}
and for $i \in \{1,\ldots, N-1 \}$, 
\begin{eqnarray}\label{ker:j<N}
\text{ker}(\rho_{i_{|_V}}) \supseteq  \big( \bigoplus\limits_{\substack{j = 1 \\ i \neq j}}^{n-1} (V \cap H^0(E_i(-p_j))) \big) \oplus (V \cap H^0(E_n(-p_1-\ldots -p_{n-1}))).
\end{eqnarray}
\end{remark}

\begin{lemma}\label{lemma for each i, kernel non zero}
Let $C$ be a comb-like curve. Let $V \cap \mathrm{H}^0(E_i(-p_i)) \neq \{ 0 \}$ and $V \cap \mathrm{H}^0(E_j(-p_j)) \neq \{ 0 \}$ for some distinct pair $i, j$  in $\{ 1, \dots , n-1 \}$. Then for all $i$, $\ker({\rho_i}_{\vert_V}) \neq 0$.
\begin{proof}
we conclude the result from \eqref{ker:N} and \eqref{ker:j<N}.
\end{proof}
\end{lemma}

\begin{example}
Let $C$ be a comb-like curve with the smooth components $C_1, \dots, C_n$.
For each $i$, let $E_i$ be a vector bundle on $C_i$ with sufficiently large degree such that $E_i$ is generated by global sections. Consider the vector bundle $E$ on $C$ obtained by glueing $E_i$ for all $i$ (one finds the glueing of vector bundles over a reducible nodal curves technique in \cite{nagaraj-seshadri96}). Then therefore $\mathrm{H}^0(C, E)$ generates $E$.

Now, take the linear subspace $V$ is the whole space $\mathrm{H}^0(C, E)$. Then for any pair $1 \le i < j \le n-1$, 
\[
\mathrm{H}^0(C, E) \cap \mathrm{H}^0(E_i(-p_i)) = \mathrm{H}^0(E_i(-p_i)) \neq 0 
\]
and
\[
\mathrm{H}^0(C, E) \cap \mathrm{H}^0(E_j(-p_j)) = \mathrm{H}^0(E_i(-p_i)) \neq 0. 
\]
By using \eqref{ker:N} and \eqref{ker:j<N}, we conclude that for each $i$, $\ker({\rho_i}_{\vert_{V}}) \neq 0$.
\end{example}
\begin{theorem}\label{thm on strongly unstabe for any rank}
Let $C$ be a comb-like curve. Let $(E, V)$ be a generated pair on $C$ with $\frac{d}{l-r} > n-1$. Suppose that $\ker({\rho_i}_{\vert_{V}}) \neq 0$, for each $i$. Then $M_{E, V}$ is strongly unstable.
 
 \begin{proof}
 Suppose on the contrary that there exists a polarization $w = (w_1, \dots , w_n)$ such that $M_{E, V}$ is $w$-semistable.
\noindent
 By using \eqref{exact sequence defining I_C_i}, we first note that  
 \begin{eqnarray*}
 \mathcal{O}_{C_i}(- p_i) \subseteq \mathcal{O}_C \quad \quad(\text{for all}\; 1 \leq i \leq n-1)
 \end{eqnarray*}
 and 
 \begin{eqnarray*}
 \mathcal{O}_{C_n} (-p_1 - \cdots - p_{n-1}) \subseteq \mathcal{O}_C.
 \end{eqnarray*}
 After tensoring with $M_{E, V}$, this  implies that 
 \begin{eqnarray}\label{eq: 1}
  {M_{E,V}}_{\vert_{C_i}}(- p_i) = M_{E, V} \otimes \mathcal{O}_{C_i}(- p_i) \subseteq M_{E, V} \otimes  \mathcal{O}_C =M_{E, V},
 \end{eqnarray}
 $\text{for all}\; 1 \leq i \leq n-1$, and 
 \begin{eqnarray*}
 {M_{E,V}}_{\vert_{C_n}}(-p_1 - \cdots - p_{n-1}) = M_{E, V} \otimes \mathcal{O}_{C_n} (-p_1 - \cdots - p_{n-1}) \subseteq M_{E, V} \otimes  \mathcal{O}_C = M_{E, V}.
 \end{eqnarray*}
 Since $\ker({\rho_i}_{\vert_{V}}) \neq 0$, from the commutative diagram \ref{diagram 1}, we have 
 \begin{eqnarray}\label{eq:09}
   \ker({\rho_i}_{\vert_{V}}) \otimes \mathcal{O}_{C_{i}} \subseteq {M_{E,V}}_{\vert_{C_i}},
\end{eqnarray}
for all $i = 1, \dots, n-1$.

\noindent
Tersoring in \eqref{eq:09} with $\mathcal{O}_{C_i}(- p_i)$, for $1 \leq i \leq n-1$, and then using \eqref{eq: 1},  we have
\begin{eqnarray*}
 \ker({\rho_i}_{\vert_{V}}) \otimes \mathcal{O}_{C_{i}}(-p_i) \subseteq {M_{E,V}}_{\vert_{C_i}}(- p_i) \subset M_{E,V} \otimes \mathcal{O}_C = M_{E, V}.
\end{eqnarray*}
Therefore, for $i = 1, \dots , n-1$, $ \ker({\rho_i}_{\vert_{V}}) \otimes \mathcal{O}_{C_{i}}(-p_i)$ is a subsheaf of $M_{E, V}$. Also, note that
\begin{eqnarray}\label{eq:10}
\begin{split}
\mu_w (\ker({\rho_i}_{\vert_{V}}) \otimes \mathcal{O}_{C_{i}}(p_i))  = \frac{t_i \chi(\mathcal{O}_{C_i}(-p_i))}{t_i w_i} = \frac{-g_i}{w_i},
\end{split}
\end{eqnarray} 
where we denote $\mathrm{rk}(\ker({\rho_i}_{\vert_{V}}))$ by $t_i$.
Now, by using $w$-semistability of $M_{E, V}$ and \eqref{eq:10}, we have the following
\begin{eqnarray*}
\begin{split}
\frac{-g_i}{w_i} & = \mu_w (\ker({\rho_i}_{\vert_{V}}) \otimes \mathcal{O}_{C_{i}}(p_i)) \\ 
              & \leq \mu_w (M_{E, V}) \\
              & = \frac{-((l-r)(p_a(C) -1) +d)}{(l-r)}.
\end{split}
\end{eqnarray*}
From the above inequalities, we have
\begin{eqnarray}\label{ineq:03}
(l-r) g_i \geq ((l-r)(p_a(C) -1) +d)w_i, \quad \quad \text{for all} \;i = 1, \dots , n-1
\end{eqnarray} 
Summing \eqref{ineq:03} over $i = 1, \dots , n-1$, we have
\begin{eqnarray}\label{ineq:04}
(l-n) \sum\limits_{i =1}^{n-1} g_i \geq ((l-r)(p_a(C) -1) +d) \sum\limits_{i =1}^{n-1}w_i
\end{eqnarray}
Similarly, we see that $ \ker({\rho_i}_{\vert_{V}}) \otimes \mathcal{O}_{C_{i}}(-p_1 - \cdots - p_{n-1})$ is a subsheaf of $M_{E, V}$, we have the following inequality
\begin{eqnarray}\label{ineq:05}
(l-r) g_n + (l-r)(n -2) \geq ((l-r)(p_a(C) -1) +d)w_n
\end{eqnarray} 
Summing \eqref{ineq:04} and \eqref{ineq:05}, we have
\begin{eqnarray*}
\begin{split}
(l-r) \sum\limits_{i =1}^{n} g_i + (l-r)(n -2) & \geq ((l-r)(p_a(C) -1) +d) \sum\limits_{i =1}^{n}w_i \\ 
& \geq (l-r)(p_a(C) -1) +d \quad \quad (\text{for} \; \sum\limits_{i =1}^{n}w_i =1).              
\end{split}
\end{eqnarray*}
Thus, we now have
\begin{eqnarray*}
(l-r) \sum\limits_{i =1}^{n} g_i + (l-r)(n -2) \ge (l-r)(p_a(C) -1) +d.
\end{eqnarray*}
By using \eqref{exact seq calculation of arithmetic genus}, we arrive at
$
n-1 \geq \frac{d}{l-r}$. This yields a contradiction of the hypothesis $n-1 < \frac{d}{l-r}$. Therefore, the bundle $M_{E, V}$ is strongly unstable. 

\end{proof}  
 \end{theorem}
 
\begin{corollary}
 Let $(E, V)$ be a generated pair on a comb-like curve $C$ with $\frac{d}{l-r} > n-1$. Let $V \cap \mathrm{H}^0(E_i(-p_i)) \neq \{ 0 \}$ and $V \cap \mathrm{H}^0(E_j(-p_j)) \neq \{ 0 \}$ for $i, \; j \in \{ 1, \dots , n-1 \}$ with $i \neq j$. Then $M_{E, V}$ is strongly unstable. 
\end{corollary}

	
\end{document}